\documentstyle[12pt]{article}
\begin{document}
\parindent 16pt
\title{\Large\bf A revised pre-order principle and  set-valued Ekeland variational  principle$^1$}
\setcounter{footnote}{1}
 \footnotetext{ This work was supported by
the National Natural Science Foundation of China (10871141).

 E-mail address: qjhsd@sina.com,  jhqiu@suda.edu.cn}

\author{  {Jing-Hui Qiu
}\\
{\footnotesize\sl School of Mathematical Sciences,  Soochow
University, Suzhou 215006, P. R. China} \\
}
\date{}
\maketitle
\begin{center}
\begin{minipage}{124mm}
\vskip 0.5cm {{\bf Abstract.}  \  In my former paper ``A pre-order
principle and set-valued Ekeland variational principle" (see: arXiv:
1311.4951[math.FA]), we established a general pre-order principle.
From the pre-order principle, we deduced most of the known
set-valued Ekeland variational principles (denoted by EVPs) and
their improvements. But the pre-order principle could not imply
Khanh and Quy's EVP in [On generalized Ekeland's variational
principle and equivalent formulations for set-valued mappings, J.
Glob. Optim., 49 (2011), 381-396], where the perturbation contains a
weak $\tau$-function. In this paper, we give a revised version of
the pre-order principle. This revised version not only implies the
original pre-order principle, but also can be applied to obtain the
above Khanh and Quy's EVP. Thus, the revised pre-order principle
implies all the known set-valued EVPs in set containing forms (to my
knowledge).\\

}\\


\end{minipage}
\end{center}
\vskip 1cm \baselineskip 18pt

\section*{ \large\bf 1. Introduction  }

\hspace*{\parindent}  In [20], we established a general pre-order
principle, which consists of a pre-order set  $(X, \preceq)$  and a
monotone extended real-valued function $\eta$ on $(X, \preceq)$. The
pre-order principle states that there exists a strong minimal point
dominated by any given point provided that the monotone function
$\eta$ satisfies three general conditions. The key to the proof of
the general pre-order principle is to distinguish two different
points by scalarizations. From the pre-order principle we obtained a
very general set-valued Ekeland variational principle (briefly,
denoted by EVP), which implies most of the known set-valued EVPs and
their improvements, for example, Ha's EVP in [7], Qiu's EVP in [17],
Bednarczuk and Zagrodny's EVP in [3], Guti\'{e}rrez, Jim\'{e}nez and
Novo's EVPs in [6],  Tammer and Z$\breve{a}$linescu's EVPs in [22],
Flores-Baz\'{a}n, Guti\'{e}rrez and  Novo's EVPs in [4], Liu and
Ng's EVPs in [13], Qiu's EVPs in [18], Khanh and Quy's EVPs in [11]
and Bao and Mordukhovich's EVPs in [1, 2]. However, it could not
imply Khanh and Quy's EVPs in [10], where the perturbations contain
 weak $\tau$-functions. The crux of the problem is that weak
$\tau$-function $p(x, x^{\prime})$ may be zero even though $x\not=
x^{\prime}$. This inspired us to find a revised version of the
pre-order principle so that it can imply Khanh and Quy's EVPs in
[10] as well.

In this paper, we  give a revised version of the general pre-order
principle. The revised version implies the original pre-order
principle in [20] and can be applied to the case that the
perturbation contains a weak $\tau$-function. Hence, the revised
pre-order principle not only implies all the set-valued EVPs in
[20], but also implies Khanh and Quy's EVPs in [10] and their
improvements. Thus, the revised pre-order principle implies all the
known
set-valued EVPs in set containing forms (to my knowledge).\\

\section*{ \large\bf 2.  A revised pre-order  principle }

\hspace*{\parindent} Let $X$ be a nonempty set. As in [4], a binary
relation $\preceq$ on $X$ is called a pre-order if it satisfies the
transitive property; a quasi order if it satisfies the reflexive and
transitive properties; a partial order if it satisfies the
antisymmetric, reflexive and transitive properties. Let $(X,
\preceq)$ be a pre-order set. An extended real-valued function
$\eta:\, (X, \preceq) \rightarrow R\cup\{\pm \infty\}$ is called
monotone with respect to $\preceq$ if for any $x_1,\, x_2\in X$,
$$x_1\preceq x_2\ \ \Longrightarrow\ \ \eta(x_1)\leq \eta(x_2).$$
For any given $x_0\in X$, denote $S(x_0)$ the set $\{x\in X:\,
x\preceq x_0\}$. We give a revised version of [20, Theorem 2.1] as follows.\\

{\bf Theorem 2.1.} \ {\sl Let $(X, \preceq)$ be a pre-order set,
$x_0\in X$ such that $S(x_0)\not=\emptyset$ and $\eta:\, (X,
\preceq) \rightarrow R\cup\{\pm\infty\}$ be an extended real-valued
function which is monotone with respect to $\preceq$.

Suppose that the following conditions are satisfied:

{\rm (A)}  \  $-\infty < \inf\{\eta(x):\, x\in S(x_0)\} <+\infty$.

{\rm (B)}  \  For any $x\in S(x_0)$ with $-\infty<\eta(x)<+\infty$
and any $z_1, z_2\in S(x)$ with $z_1\not= z_2$, one has $\eta(x)>
\min\{\eta(z_1),\eta(z_2)\}$.

{\rm (C)} \ For any sequence $(x_n)\subset S(x_0)$ with $x_n\in
S(x_{n-1}),\ \forall n$, such that  $\eta(x_n) -\inf\{\eta(x):\,
x\in S(x_{n-1})\} \rightarrow 0$ \  $(n\rightarrow \infty)$, there
exists $u\in X$ such that $u\in S(x_n),\ \forall n$.

Then there exists $\hat{x}\in X$ such that

{\rm (a)} \ $\hat{x}\in S(x_0)$;

{\rm (b)} \ $S(\hat{x})\subset \{\hat{x}\}$.}\\

{\bf Proof.} \  The first half of the proof is similar to that of
[20 , Theorem 2.1]. Here, for the sake of completeness, we give the
whole process.

For brevity, we denote $\inf\{\eta(x):\, x\in S(x_0)\}$ by $\inf
\eta\circ S(X_0)$. By (A), we have
$$-\infty <\inf \eta\circ S(x_0)<+\infty. \eqno{(2.1)}$$
So, there exists $x_1\in S(x_0)$ such that
$$\eta(x_1) < \inf \eta\circ S(x_0) +\frac{1}{2}. \eqno{(2.2)}$$
By the transitive property of $\preceq$, we have
$$ S(x_1)\subset S(x_0). \eqno{(2.3)}$$
If $S(x_1)\subset \{x_1\}$, then we may take $\hat{x}:=x_1$ and
clearly $\hat{x}$ satisfies (a) and (b). If not, by (2.1), (2.2) and
(2.3) we conclude that
$$-\infty < \inf \eta\circ S(x_1) <+\infty.$$
So, there exists $x_2\in S(x_1)$ such that
$$\eta(x_2) < \inf\eta\circ S(x_1) +\frac{1}{2^2}.$$
In general, if $x_{n-1}\in X$ has been chosen, we may choose $x_n\in
S(x_{n-1})$ such that
$$\eta(x_n) < \inf\eta\circ S(x_{n-1}) +\frac{1}{2^n}.$$
If there exists $n$ such that $S(x_n)\subset \{x_n\}$, then we may
take $\hat{x} :=x_n$ and clearly $\hat{x}$ satisfies (a) and (b). If
not, we can obtain a sequence $(x_n)\subset S(x_0)$ with $x_n\in
S(x_{n-1}),\ \forall n$, such that
$$\eta(x_n) < \inf \eta\circ S(x_{n-1}) +\frac{1}{2^n},\ \forall n.
\eqno{(2.4)}$$ Obviously, $\eta(x_n) -\inf\eta\circ S(x_{n-1})
\rightarrow 0$ when $n\rightarrow \infty$. By (C), there exists
$u\in X$ such that $$u\in S(x_n), \ \ \forall n. \eqno{(2.5)}$$
Obviously, $u\in S(x_0)$, i.e., $u$ satisfies (a).

By (A), $\eta(u)\geq \inf\,\eta\circ S(x_0) >-\infty.$ Also, from
$u\in S(x_1)$ and (2.2), we have
$$\eta(u) \leq \eta(x_1) < \inf\,\eta\circ S(x_0) + \frac{1}{2}
<+\infty.$$ Hence $$-\infty < \eta(u) <+\infty.\eqno{(2.6)}$$ Now,
we assert that $S(u)$ could not contain two different points. If
not, there exists $z_1, z_2\in S(u)$ such that $z_1\not=z_2$. By
(2.6) and (B), we have
$$\eta(u) > \min\{\eta(z_1), \eta(z_2)\}.$$
For definiteness, assume that
$$\eta(u) > \eta(z_1).\eqno{(2.7)}$$
On the other hand, by  $u\in S(x_n)$ and $z_1\in S(u)$, we have
$$z_1\in S(u) \subset S(x_n), \ \ \forall n.\eqno{(2.8)}$$
By $u\preceq x_n$, (2.4) and (2.8), we have
$$\eta(u) \leq \eta(x_n) < \inf\,\eta\circ S(x_{n-1})
+\frac{1}{2^n}\leq \eta(z_1) +\frac{1}{2^n}.$$ Letting
$n\rightarrow\infty$, we have $\eta(u)\leq \eta(z_1)$, which
contradicts (2.7). Thus, we have shown that $S(u)$ is only  empty
set or singleton.

If $S(u)=\emptyset$, then put $\hat{x}:=u$. Obviously, $\hat{x}$
satisfies (a) and (b).

If $S(u)$ is a singleton $\{v\}$, then we conclude that $S(v)\subset
\{v\}$. In fact, if there exists $w\in S(v)\backslash\{v\}$, then
$w\not=v$ and $v, w\in S(u)$. This contradicts the conclusion that
$S(u)$ could not contain two different points. Also, we remark that
$v\in S(u) \subset S(x_0)$. Thus, $\hat{x}:=v$ satisfies (a) and
(b).\hfill\framebox[2mm]{}\\

{\bf Corollary  2.1.} \ (see [20, Theorem 2.1]) \ {\sl The result of
Theorem 2.1 remains true if assumption {\rm (B)} is replaced by the
following assumption:

{\rm (${\rm B_0}$)} \ For any $x\in S(x_0)$ with $-\infty <\eta(x)
<+\infty$ and any $x^{\prime}\in S(x)\backslash\{x\}$, one has
$\eta(x) >\eta(x^{\prime})$.}\\

{\bf Proof.} \  It is sufficient to show that (${\rm B_0}$)\
$\Rightarrow$\ (B). Let $x\in S(x_0)$ with $-\infty <\eta(x)
<+\infty$ and let $z_1, z_2 \in S(x)$ with $z_1\not= z_2$.
Obviously, at least one of $z_1$ and $z_2$ is not equal to $x$. For
definiteness, we assume that $z_1\not= x$. Thus, by (${\rm B_0}$) we
have $\eta(x) > \eta(z_1)$. Certainly, $\eta(x) > \min\{\eta(z_1),
\eta(z_2)\}$ and (B) holds.\hfill\framebox[2mm]{}\\

{\bf Remark 2.1.} \  From Corollary 2.1, we know that Theorem 2.1
implies [20, Theorem 2.1]. Hence, from Theorem 2.1 we can deduce all
the set-valued EVPs in [20]. Furthermore, we shall see that Theorem
2.1 also implies the set-valued EVPs in [10] and their
improvements.\\

\section*{ \large\bf 3.  A general  set-valued EVP  }

\hspace*{\parindent} Let  $Y$ be a real linear space. If $A,\,
B\subset Y$ and $\alpha\in R$, the sets $A+B$ and $\alpha\,A$ are
defined as follows: $$A+B:=\{z\in Y:\, \exists x\in A,\, \exists
y\in B\ {\rm such\ that}\  z=x+y\},$$  $$\alpha A:=\{z\in Y:\,
\exists x\in A\ \ {\rm such\ that}\ \  z=\alpha x\}.$$ A nonempty
subset $D$ of $Y$ is called a cone if $\alpha D\subset D$ for any
$\alpha\geq 0$. And $D$ is called a convex cone if $D+D\subset D$
and $\alpha D\subset D$ for any $\alpha\geq 0$. A convex cone $D$
can specify a quasi order $\leq_D$ on $Y$ as follows:
$$y_1,\, y_2\in Y,\ \ y_1\leq_D y_2\ \ \ \Longleftrightarrow\ \ \ \ y_1-y_2\in
-D.$$ In this case, $D$ is also called the ordering cone or positive
cone. We always assume that $D$ is nontrivial, i.e., $D\not=\{0\}$
and $D\not=Y$. An extended real function $\xi:\, Y\rightarrow
R\cup\{\pm\infty\}$ is said to be $D$-monotone if
$\xi(y_1)\leq\xi(y_2)$ whenever $y_1\leq_D y_2$. For any nonempty
subset $M$ of $Y$, we put $\inf \xi\circ M\,=\,\inf\{\xi(y):\, y\in
M\}$. If $\inf \xi\circ M
>-\infty$, we say that $\xi$ is lower bounded on $M$. For any given
$y\in Y$, sometimes we denote $\xi(y)$ by $\xi\circ y$. A family of
set-valued maps $F_{\lambda}:\, X\times X\rightarrow
2^D\backslash\{\emptyset\},\ \lambda\in\Lambda,$ is said to satisfy
the ``triangle inequality" property (briefly, denoted by property
TI, see [4]) if for each $x_i\in X,\ i=1,2,3,$ and
$\lambda\in\Lambda$ there exist $\mu,\nu\in\Lambda$ such that
$$F_{\mu}(x_1, x_2) + F_{\nu}(x_2, x_3)\,\subset\, F_{\lambda}(x_1,
x_3) +D.$$ Let $X$ be a nonempty set and let $f:\, X\rightarrow
2^Y\backslash\{\emptyset\}$ be a set-valued map. For any nonempty
set $A\subset X$, we put $f(A):=\cup\{f(x):\, x\in A\}$. For any
$x_1,\, x_2\in X$, define $x_2\preceq x_1$ iff
$$f(x_1)\subset f(x_2) + F_{\lambda}(x_2, x_1) + D,\ \
\forall\lambda\in\Lambda.$$\\

{\bf Lemma 3.1.}\ (see [20, lemma 3.1]) \ {\sl ``$\preceq$" is a
pre-order on $X$, i.e., it
is a binary relation satisfying transitive property.}\\

Now, we can prove the following theorem, which is indeed a
generalization of [20, Theorem 3.1].\\

{\bf Theorem 3.1.} \ {\sl Let $X$ be a nonempty set, $Y$ be a real
linear space, $D\subset Y$ be a convex cone specifying a quasi order
$\leq_D$ on $Y$, $f:\, X\rightarrow 2^Y\backslash\{\emptyset\}$ be a
set-valued map and $F_{\lambda}:\, X\times X\rightarrow
2^D\backslash\{\emptyset\},\ \lambda\in \Lambda$, be a family of
set-valued maps satisfying the property TI. Let $x_0\in X$ such that
$$S(x_0):=\{x\in X:\, f(x_0)\subset f(x) + F_{\lambda}(x, x_0) +D,\
\forall\lambda\in\Lambda\}\,\not=\,\emptyset.$$ Suppose that there
exists a $D$-monotone extended real-valued function $\xi:\,
Y\rightarrow R\cup\{\pm\infty\}$ satisfying the following
assumptions:

{\rm (D)} \  $-\infty  <\inf \xi\circ f(S(x_0))<+\infty$.

{\rm (E)} \   For any $x\in S(x_0)$ with $-\infty< \inf \xi\circ
f(x)<+\infty$ and for any $z_1, z_2\in S(x)$ with $z_1\not= z_2$,
one has $\inf\xi\circ f(x)
> \min\{\inf \xi\circ f(z_1), \inf \xi\circ f(z_2)\}$.

{\rm (F)} \   For any  sequence $(x_n)\subset S(x_0)$ with $x_n\in
S(x_{n-1}),\ \forall n,$ such that  $\inf \xi\circ f(x_n)
-\inf\xi\circ f(S(x_{n-1})) \rightarrow 0\ \,(n\rightarrow\infty)$,
there exists $u\in X$ such that $u\in S(x_n),\ \forall n.$

Then there exists $\hat{x}\in X$ such that

{\rm (a)} $f(x_0)\subset f(\hat{x}) +F_{\lambda}(\hat{x}, x_0) +D,\
\forall\lambda\in\Lambda$;

{\rm (b)} $\forall x\in X\backslash\{\hat{x}\}, \,
\exists\lambda\in\Lambda$ \ such \ that\ $f(\hat{x})\not\subset f(x)
+F_{\lambda}(x,\hat{x}) +D.$ }\\

{\bf Proof.}\  By Lemma 3.1, we can define a pre-order $\preceq$ on
$X$ as follows: for any $x_1,\, x_2\in X$,
$$x_2\preceq x_1 \ \ \  {\rm iff}\ \ \  f(x_1)\subset f(x_2) +
F_{\lambda}(x_2, x_1) +D,\ \forall\lambda\in\Lambda.$$ Thus,
$S(x_0)=\{x\in X:\, x\preceq x_0\}$. Define an extended real-valued
function  $\eta:\, (X, \preceq) \rightarrow R\cup\{\pm\infty\}$ as
follows
$$\eta(x):=\inf \xi\circ f(x),\ \ \forall x\in X.$$
Let $x^{\prime}\preceq x$. Then
$$f(x) \subset f(x^{\prime}) + F_{\lambda}(x^{\prime}, x) +D,\ \
\forall\lambda\in\Lambda.$$ For any $y\in f(x)$, there exists
$y^{\prime}\in f(x^{\prime}),\ q_{\lambda}(x^{\prime}, x) \in
F_{\lambda}(x^{\prime}, x),\ d_{\lambda, x^{\prime}, x}\in D$ such
that
$$y= y^{\prime} + q_{\lambda}(x^{\prime}, x) + d_{\lambda,
x^{\prime}, x} .$$ Since $$y-y^{\prime} = q_{\lambda}(x^{\prime}, x)
+d_{\lambda, x^{\prime}, x} \in D,\ {\rm i.e.,}\ y\geq_D
y^{\prime},$$ we have
$$\xi(y) \geq \xi(y^{\prime}) \geq \inf\xi\circ f(x^{\prime}).$$
As $y\in f(x)$ is arbitrary, we have
$$\inf \xi\circ f(x) \geq \inf \xi\circ f(x^{\prime}),\ \  {\rm
i.e.,}\  \  \eta(x)\geq\eta(x^{\prime}).$$ Thus, $\eta$ is monotone
with respect to $\preceq$. It is easy to see that assumptions (D),
(E) and (F) are exactly assumptions (A), (B) and (C) in Theorem 2.1.
Now, applying Theorem 2.1, we know that there exists $\hat{x}\in X$
such that $\hat{x}\in S(x_0)$ and $S(\hat{x})\subset \{\hat{x}\}$.
This means that
$$f(x_0)\subset f(\hat{x}) +F_{\lambda}(\hat{x}, x_0)
+D,\ \forall\lambda\in\Lambda$$ and
 $$\forall x\in X\backslash\{\hat{x}\}, \, \exists\lambda\in\Lambda
\ {\rm such \ that}\ f(\hat{x})\not\subset f(x)
+F_{\lambda}(x,\hat{x}) +D.$$ That is, $\hat{x}$ satisfies (a) and (b). \hfill\framebox[2mm]{}\\

For a real linear space $Y$, denote the algebraic dual of $Y$ by
$Y^{\#}$ and  denote the positive polar cone of  $D$ in $Y^{\#}$ by
$D^{+\#}$, i.e., $D^{+\#}=\{l\in Y^{\#}:\, l(d)\geq 0,\,\forall d\in
D\}$. Particularly, if the $\xi$ in Theorem 3.1 is an element of
$D^{+\#}\backslash\{0\}$. then assumptions (D) and (E) become more
concise.\\

{\bf Theorem  3.1$^{\prime}$.} \ {\sl  let $X,\,Y,\, D,\, f,\,
F_{\lambda},\,\lambda\in\Lambda$, and $x_0\in X$ be the same as in
Theorem 3.1. Suppose that there exists $\xi\in
D^{+\#}\backslash\{0\}$ satisfying the following assumptions:

{\rm (D)} \  $\xi$ is lower bounded on $f(S(x_0))$, i.e., $-\infty <
\inf \xi\circ f(S(x_0))$.

{\rm (E)} \  For any $x\in S(x_0)$ and any  $z_1, z_2\in S(x)$ with
$z_1\not= z_2$, one has  $\inf \xi\circ f(x) > \min \{\inf \xi\circ
f(z_1), \inf\xi\circ f(z_2)\}$.

{\rm (F)} \ See Theorem 3.1.

Then the result of Theorem 3.1 remains true.}\\

Obviously, Theorem 3.1$^{\prime}$ here is a generalization of [20,
Theorem 3.1$^{\prime}$]. Hence, from Theorem 3.1$^{\prime}$ we can
deduce a number of specific versions of EVP in [20]. Moreover, it is
encouraging that Theorem 3.1$^{\prime}$ can be applied to the case
that the perturbation contains a weak $\tau$-function. From this, we
can deduce Khanh and Quy's EVPs in [10] and their generalizations.\\

{\section*{ \large\bf 4.   Set-valued EVPs where perturbations
contain p-distance, q-distances and r-distances}

\hspace*{\parindent}  In order to derive Khanh and Quy's EVPs in
[10], we need recall the definitions of $\tau$-functions and weak
$\tau$-functions.\\

{\bf Definition 4.1.} \ (see [10, 12]) \ Let $(X, d)$ be a metric
space.  Then a real-valued function $p:\, X\times X \rightarrow [0,
+\infty)$ is called a $\tau$-function on $X$ if the following
conditions are satisfied:

($\tau$1) \  for any $x, y, z\in X,\ p(x,z)\leq p(x, y) +p(y,z)$;

($\tau$2) \  if $x\in X$ and a sequence $(y_n)$ with $y_n\rightarrow
y$ in $X$ and $p(x, y_n)\leq M$ for some $M=M(x) >0$, then $p(x,y)
\leq M$;

($\tau$3) \ for any sequence $(y_n)$ in $X$ with $p(y_n,
y_m)\rightarrow 0\ (m>n\rightarrow\infty)$, if there exists $(z_n)$
in $X$ such that $\lim\limits_{n\rightarrow\infty} p(y_n ,z_n) =0$,
then $\lim\limits_{n\rightarrow\infty} d(y_n, z_n) =0$;

($\tau$4) \ for $x,y,z\in X,\ p(z,x) =0$ and $p(z,y) =0$ imply
$x=y$.

A real-valued function $p:\, X\times X \rightarrow [0, +\infty)$ is
called a weak $\tau$-function if condition ($\tau$2) is removed.
i.e., only conditions ($\tau$1), ($\tau$3) and ($\tau$4) hold.\\

It is known ([12]) that a w-distance (see [9]) is a $\tau$-function.
In [21], we introduced the notions of p-distances and q-distances.
Here, we further introduce a more general notion: r-distances. We
list these notions as follows.\\

{\bf Definition 4.2.} \  Let $X$  be a separated uniform space. An
extended real-valued function $p:\, X\times X \rightarrow [0,
+\infty]$ is called a p-distance on $X$ if the following conditions
(q1)--(q3) are satisfied:

(q1) \ for any $x, y, z\in X,\ p(x,z)\leq p(x,y) +p(y,z)$;

(q2) \ every sequence $(y_n)\subset X$ with $p(y_n, y_m) \rightarrow
0\ (m>n\rightarrow \infty)$ is a Cauchy sequence and in the case
$p(y_n, y) \rightarrow 0$ is equivalent to $y_n\rightarrow y$ in
$X$;

(q3) \ for $x,y, z\in X,\ p(z,x) =0$ and $p(z,y) =0$ imply $x=y$.

If condition (q2) is replaced by the following weaker condition

(q2$^{\prime}$) \ every sequence $(y_n)\subset X$ with $p(y_n,
y_m)\rightarrow 0\ (m>n\rightarrow\infty)$ is a Cauchy sequence and
in the case $p(y_n, y) \rightarrow 0$ implies $y_n\rightarrow y$ in
$X$,

then $p$ is called a q-distance on $X$.

Moreover, if condition (q2) is replaced by the following further
weaker condition

(q2$^{\prime\prime}$) \ every sequence $(y_n)\subset X$ with $p(y_n,
y_m) \rightarrow 0\ (m>n\rightarrow\infty)$ is a Cauchy sequence,

then $p$ is called an r-distance on $X$.\\

Obviously, every weak $\tau$-function is a finite-valued q-distance.
Certainly, it is also an r-distance.\\

{\bf Definition 4.3.} \ (see [18, Definition 3.4])  \ Let $X$ be a
topological space and let $S(\cdot): X\rightarrow
2^X\backslash\{\emptyset\}$ be a set-valued map. The set-valued map
$S(\cdot)$ is said to be dynamically closed at $x\in X$ if
$(x_n)\subset S(x),\ S(x_{n+1})\subset S(x_n)\subset S(x)$ for all
$n$ and $x_n\rightarrow \bar{x}$, then $\bar{x}\in S(x)$. In this
case, we also say that $S(x)$ is dynamically closed. Moreover, let
$X$ be an uniform space and $x\in X$. Then $X$ is said to be
$S(x)$-dynamically complete if every Cauchy
sequence $(x_n)\subset S(x)$ such that $S(x_{n+1})\subset S(x_n)\subset S(x)$ for all $n$, is convergent in $X$.\\

If $Y$ is a separated topological vector space, we denote $Y^*$ the
topological dual of $Y$ and denote $D^+$ the positive polar cone of
$D$ in $Y^*$, i.e., $D^+:=\{l\in Y^*:\, l(d)\geq 0,\ \forall d\in
D\}$. It is possible that $Y^*=\{0\}$, i.e., there is no non-trivial
continuous linear functional on $Y$. However, if $Y$ is a locally
convex separated topological vector space (briefly, denoted by a
locally convex space), then $Y^*$ is large enough so that $Y^*$ can
separate points in $Y$. for a convex set $H\subset
D\backslash\{0\}$, we  denote $H^{+s}$ the set $\{l\in Y^*:\,
\inf\,l\circ H >0\}$. It is easy to see that $H^{+s}\cap D^+
\not=\emptyset$ iff $0\not\in {\rm cl}(H+D)$. Now, we can give our
first special set-valued EVP, where the perturbation consists of an
r-distance and a convex subset $H$ of the ordering cone $D$.\\

{\bf Theorem 4.1.}   \ {\sl Let $(X, {\cal U})$ be a sequentially
complete uniform space,  $Y$ be a locally convex space quasi ordered
by a convex cone $D$, $H\subset D\backslash\{0\}$ be a convex set,
$f:\, X\rightarrow 2^Y\backslash \{\emptyset\}$ be a set-valued map
and $p: X\times X \rightarrow [0,+\infty)$ be an r-distance. Let
$x_0\in X$ such that $S(x_0):=\{x\in X:\, f(x_0)\subset f(x) +
p(x_0, x) H +D\} \not= \emptyset$ and for any $x\in S(x_0)$, $S(x)$
is dynamically closed.

Suppose that the following assumptions are satisfied:

{\rm (B1)} \ $H^{+s}\cap D^+ \not= \emptyset$, or equivalently,
$0\not\in {\rm cl}(H+D)$.

{\rm (B2)} \  $\exists\, \xi\in H^{+s}\cap D^+$ such that $\xi$ is
lower bounded on $f(S(x_0))$.

Then, there exists $\hat{x}\in X$ such that

{\rm (a)}  $f(x_0)\subset f(\hat{x}) +p(x_0, \hat{x}) H +D$;

{\rm (b)} $\forall x\in X\backslash\{\hat{x}\},
f(\hat{x})\not\subset f(x) + p(\hat{x}, x) H +D.$

If $p(x, x)= 0$ for all $x\in X$, then $\hat{x}$ is a strict
minimizer of $f(\cdot)  + p(\hat{x}, \cdot) H$. {\rm (concerning strict minimizer, see [7, 10])}}\\

{\bf Proof.} \  Put $F(x, x^{\prime}) = p(x^{\prime}, x) H,\ \forall
x, x^{\prime}\in X$. Clearly, the family $\{F\}$ satisfies property
TI. By (B2), there exists $\xi\in H^{+s}\cap D^+$ such that $\xi$ is
lower bounded on $f(S(x_0))$. That is, assumption (D) in Theorem
3.1$^{\prime}$ is satisfied. Denote $\inf\{\xi(h):\, h\in H\}$ by
$\alpha$. Then $\alpha >0$. We have: $\xi(h)\geq \alpha,\ \forall
h\in H$ and $\xi(d)\geq 0,\ \forall d\in D.$ Let $x\in S(x_0)$ and
let $z_1, z_2\in S(x)$ with $z_1\not= z_2$. Then
$$f(x) \subset  f(z_i) + p(x, z_i)H +D,\ \ i=1,2.\eqno{(4.1)}$$
By the definition of r-distances, at least one of $p(x, z_1)$ and
$p(x, z_2)$ is strictly greater than $0$. For definiteness, we
assume that $p(x, z_1) >0$. From (4.1), we easily see that
$$\inf\,\xi\circ f(x) \geq \inf\,\xi\circ f(z_1) + p(x, z_1)\alpha
>\inf\,\xi\circ f(z_1).$$
Hence, $$\inf\,\xi\circ f(x) >  \min\{\inf\,\xi\circ
f(z_1),\inf\,\xi\circ f(z_2)\}.$$ That is, assumption (E) in Theorem
3.1$^{\prime}$ is satisfied.

Next, we show that assumption (F) in Theorem 3.1$^{\prime}$ is
satisfied. Let a sequence $(x_n)\subset S(x_0)$ such that $x_n\in
S(x_{n-1})$ and
$$\inf\,\xi\circ f(x_n)  < \inf\,\xi\circ f(S(x_{n-1})) +\epsilon_n\
\ \forall n,$$ where every $\epsilon_n >0$ and
$\epsilon_n\rightarrow 0$. For each $n$, take $y_n\in f(x_n)$ such
that
$$\xi\circ y_n < \inf\,\xi\circ f(S(x_{n-1})) + \epsilon_n.
\eqno{(4.2)}$$ When $m>n$, $x_m\in S(x_{m-1})\subset S(x_n)$. Hence,
$$y_n\in f(x_n)\subset f(x_m) +p(x_n, x_m) H +D.$$
Thus, there exists $y_{m,n}\in f(x_m),\ h_{m,n}\in H$ and
$d_{m,n}\in D$ such that
$$y_n = y_{m,n} + p(x_n, x_m) h_{m,n}  +d_{m,n}.$$
Acting upon two sides of the above equality by $\xi$, we have
$$\xi\circ y_n = \xi\circ y_{m,n} + p(x_n, x_m)\, \xi\circ h_{m,n}
+\xi\circ d_{m,n} \geq \xi\circ y_{m,n} + p(x_n, x_m)\,
\alpha.\eqno{(4.3)}$$ Observing that $y_{m,n}\in f(x_m) \subset
f(S(x_{n-1}))$ and using (4.3) and (4.2), we have
\begin{eqnarray*}
p(x_n, x_m) \,&\leq&\,\frac{1}{\alpha}(\xi\circ y_n -\xi\circ
y_{m,n})\\
&\leq&\, \frac{1}{\alpha}(\xi\circ y_n -\inf\xi\circ f(S(x_{n-1}))\\
&<&\, \frac{1}{\alpha}\epsilon_n.
\end{eqnarray*}
Hence, $p(x_n, x_m) \rightarrow\ 0\ \ (m>n\rightarrow\infty)$. Since
$p$ is an r-distance, $(x_n)$ is a Cauchy sequence in $X$. And since
$X$ is sequentially complete, there exists $u\in X$ such that
$x_n\rightarrow u$. For each given $n$, $x_{n+i}\rightarrow u\
(i\rightarrow\infty)$, where $(x_{n+i})_{i\in N}\subset S(x_n)$ and
$x_{n+i}\in S(x_{n+i-1}), \ \forall i$.  Remarking that $S(x_n)$ is
dynamically closed, we have $u\in S(x_n)$. Thus, assumption (F) is
satisfied. Now, applying Theorem 3.1$^{\prime}$ we obtain the
result.\hfill\framebox[2mm]{}\\

{\bf Remark 4.1} \  In particular, if the sequentially complete
uniform space $X$ is a complete metric space $(X,d)$, then we obtain
a generalization of [10, Theorem 4.1]. We see that the point $k_0$
in the ordering cone $D$ and the weak $\tau$-function $p$ appeared
in the perturbation in [10, Theorem 4.1] are replaced here by a
subset $H$ of the ordering cone $D$ and an r-distance $p$,
respectively. The condition that $k_0\not\in -{\rm cl}(D)$ is
replaced by (B1), i.e.,  $0\not\in {\rm cl}(H+D)$, or equivalently,
$H^{+s}\cap D^+ \not=\emptyset$. And the condition that $f(X)$ is
quasi-bounded from below, i.e., there exists a bounded set $M$ such
that $f(X)\subset M+D$, is replaced by (B2), i.e., there exists
$\xi\in H^{+s}\cap D^+$ such that $\xi$  is
lower bounded on $f(S(x_0))$.\\

Carefully checking the proof of Theorem 4.1, we have the following
theorem, which extends [10, Theorem 4.3].\\

{\bf Theorem 4.2.}   \ {\sl The result of Theorem 4.1 remains true
if the condition that $(X, {\cal U})$ is sequentially complete is
replaced by one that $(X, {\cal U})$ is $S(x_0)$-dynamically
complete.}\\

Traditionally, the statement of EVP, say for a scalar function $f$
on a metric space, includes an $\epsilon
>0$ such that  $f(x_0) < \inf\limits_{x\in {\rm dom}f} f(x)
+\epsilon$ and an estimate of $d(x_0,\hat{x})$. By modifying Theorem
4.2, we can obtain a corresponding statement, which extends [10,
Theorem 4.4].\\

{\bf Theorem 4.3.}   \ {\sl Let, additionally to the assumptions of
Theorem 4.2, that $f(x_0)\not\subset f(x) +\epsilon H +D$ for some
$\epsilon >0$ and all $x\in X$. Then, for any $\lambda >0$, there
exists $\hat{x}\in X$ such that

{\rm (a)} \ $f(x_0)\subset f(\hat{x}) + ({\epsilon}/{\lambda})
p(x_0, \hat{x}) H+D$;

{\rm (b)} \  $\forall x\in X\backslash\{\hat{x}\}$,
$f(\hat{x})\not\subset f(x) + (\epsilon/\lambda) p(\hat{x}, x) H
+D$;

{\rm (c)} \ $p(x_0, \hat{x}) \leq \lambda$. }\\

{\bf Proof.} \  Replacing $p$ by $(\epsilon/\lambda)p$ in the proof
of Theorem 4.1, we can obtain $\hat{x}\in X$ such that (a) and (b)
are satisfied. We claim that $p(x_0, \hat{x}) \leq\lambda$. Indeed,
if $p(x_0, \hat{x}) >\lambda$, then by (a) we would have
$$f(x_0) \subset f(\hat{x}) + \frac{\epsilon}{\lambda} p(x_0,
\hat{x}) H +D \subset f(\hat{x}) +\epsilon H +D,$$ which contradicts
the assumption that $f(x_0)\not\subset f(x) +\epsilon H +D$ for all $x\in X$.\hfill\framebox[2mm]{}\\

In order to obtain more special versions of set-valued EVP, we need
to recall some concepts. As we know, the lower semi-continuity of
scalar functions can be extended to set-valued maps. Let $X$ be a
metric space, $Y$ be a locally convex space and $D\subset Y$ be a
closed convex cone. As in [7], a set-valued map $f:\, X\rightarrow
2^Y\backslash\{\emptyset\}$ is said to be $D$-lower semi-continuous
(briefly, denoted by $D$-l.c.s.) on $X$ if for any $y\in Y$, the set
$\{x\in X:\, f(x)\cap(y-D)\not=\emptyset\}$ is closed. In [3, 5, 8,
14], sequentially lower monotone vector-valued maps were considered
(in [3], they are called monotonically semicontinuous functions; and
in [5], they are called functions satisfying (H4)). We extend the
concept to set-valued maps (see [17]). A set-valued map
$f:\,X\rightarrow 2^Y\backslash\{\emptyset\}$ is said to be
$D$-sequentially lower monotone (briefly, denoted by $D$-s.l.m.) if
$f(x_n)\subset f(x_{n+1}) +D,\ \forall n,$ and $ x_n\rightarrow
\bar{x}$ imply that $f(x_n) \subset f(\bar{x}) +D,\ \forall n$ (in
[10], such a set-valued map $f$ is said to be $D$-lower
semi-continuous from above, denoted by $D$-lsca.). Just as in the
case of vector-valued functions, we can show that  a $D$-l.s.c.
set-valued map is $D$-s.l.m., but the converse is not true (see
[17]). Here, we discuss problems in a more general setting. Let $X$
be a uniform space and $D\subset Y$ be a convex cone (needn't be
closed). In the more general setting, we still define $D$-s.l.m.
set-valued maps as above-mentioned. A set-valued map $f:\,
X\rightarrow 2^Y\backslash\{\emptyset\}$ is said to have $D$-closed
(resp., $D$-locally closed) values if for any $x\in X$, $f(x) +D$ is
closed (resp., locally closed). Concerning locally closed sets, see,
e.g., [15]. Let $H$ be a nonempty subset of a locally convex space
$Y$. A convex series of points in $H$ is a series of the form
$\sum\limits_{n=1}^{\infty} \lambda_n x_n$, where every $x_n\in H,\
\lambda_n \geq 0$ and $\sum\limits_{n=1}^{\infty}\lambda_n =1$. $H$
is said to be $\sigma$-convex if every convex series of its points
converges to a point of $H$ (see [15, Definition 2.1.4]). Sometimes,
a $\sigma$-convex set is called a cs-complete set, see [22, 23]. It
is easy to see that a $\sigma$-convex set must be bounded and
convex, but it may be non-closed. For example, an open ball in a
Banach space is a $\sigma$-convex set, but it is not closed. For
details, see, e.g., [19] and the references therein.

The following theorem is a generalization of [20, Theorem 4.2],
where the perturbation $\gamma d(x, \hat{x})H$ is replaced by a more
general version: $p(\hat{x}, x) H$. Hence, it is also a
generalization of [22, Theorem 6.2], [18, Theorem 6.8] and [13, Theorem 3.5(ii)].\\

{\bf Theorem 4.4.} \ {\sl Let $(X, {\cal U})$ be a uniform space,
$Y$ be a locally convex space quasi ordered by a convex cone $D$,
$H\subset D\backslash\{0\}$ be a convex set, $f:\, X\rightarrow
2^Y\backslash\{\emptyset\}$ be a set-valued map and $p:\, X\times
X\rightarrow [0, +\infty)$ be a p-distance. Let $x_0\in X$ such that
$S(x_0):=\{x\in X:\, f(x_0)\subset f(x) + p(x_0, x) H +D\} \not=
\emptyset$ and $(X, {\cal U})$ be $S(x_0)$-dynamically complete.
Suppose that the following assumptions are satisfied:

{\rm (B1)} \ $H^{+s}\cap D^+ \not= \emptyset$.

{\rm (B2)} \  $\exists\, \xi\in H^{+s}\cap D^+$ such that $\xi$ is
lower bounded on $f(S(x_0))$.

{\rm (B3)} \  $f$ is $D$-s.l.m.

{\rm (B4)} \   $H$ is $\sigma$-convex and $f$ has $D$-closed values
(or $D$-locally closed values).

Then, there exists $\hat{x}\in X$ such that

{\rm (a)}  $f(x_0)\subset f(\hat{x}) +p(x_0, \hat{x}) H +D$;

{\rm (b)} $\forall x\in X\backslash\{\hat{x}\},
f(\hat{x})\not\subset f(x) + p(\hat{x}, x) H +D.$}\\

{\bf Proof.} \  Put $F(x, x^{\prime}) = p(x^{\prime}, x) H,\ \forall
x, x^{\prime}\in X$. As in the proof of Theorem 4.1, by (B1) and
(B2) we can prove that assumptions (D) and (E) in Theorem
3.1$^{\prime}$ are satisfied. It is sufficient to show that
assumption (F) in Theorem 3.1$^{\prime}$ is satisfied.

Let a sequence $(x_n) \subset S(x_0)$ such that $x_n\in S(x_{n-1})$
and
$$\inf\,\xi\circ f(x_n)  < \inf\,\xi\circ f(S(x_{n-1})) +\epsilon_n\
\ \forall n,$$ where every $\epsilon_n >0$ and
$\epsilon_n\rightarrow 0$. As in the proof of Theorem 4.1, we can
show that $p(x_n, x_m) \rightarrow 0\ (m>n\rightarrow\infty).$ By
the definition of p-distances, $(x_n)$ is a Cauchy sequence in $(X,
{\cal U})$. Since $(X, {\cal U})$ is $S(x_0)$-dynamically complete,
there exists $u\in X$ such that $x_n\rightarrow u$. Again using the
definition of p-distances, we have $p(x_n, u) \rightarrow 0\
(n\rightarrow\infty)$. Take any fixed $n_0\in N$ and put $z_1:=
x_{n_0}$. As $p(x_k, u) \rightarrow 0\ (k\rightarrow\infty)$, we may
choose a sequence $(z_n)$ from $(x_k)$ such that $p(z_{n+1}, u) <
1/(n+1)$ and $z_{n+1}\in S(z_n),\ \forall n$. Take any $v_1\in
f(z_1)$. As $z_2\in S(z_1)$, we have
$$v_1\in f(z_1)\subset f(z_2) + p(z_1, z_2)H +D.$$
Hence, there exists $v_2\in f(z_2),\ h_1\in H$ and $d_1\in D$ such
that
$$v_1=v_2 + p(z_1, z_2) h_1 +d_1.$$
In general, if $v_n\in f(z_n)$ is given, then
$$v_n\in f(z_n)\subset f(z_{n+1})+ p(z_n, z_{n+1}) H +D.$$
Hence,  there exists $v_{n+1}\in f(z_{n+1}),\ h_n\in H$ and $d_n\in
D$ such that
$$v_n = v_{n+1} + p(z_n, z_{n+1}) h_n +d_n.$$
Adding two sides of the above $n$ equalities, we have
$$\sum\limits_{i=1}^n v_i\,=\,\sum\limits_{i=2}^{n+1} v_i +
\sum\limits_{i=1}^n p(z_i, z_{i+1}) h_i +\sum\limits_{i=1}^n d_i.$$
From this,
$$v_1\,=\, v_{n+1} + \sum\limits_{i=1}^n p(z_{i}, z_{i+1}) h_i
+\sum\limits_{i=1}^n d_i.\eqno{(4.4)}$$ As $\xi\in H^{+s}\cap D^+$,
$\xi(D)\geq 0$ and there exists $\alpha
>0$ such that $\xi(H)\geq\alpha$. Acting on two sides of (4.4) by
$\xi$, we have
\begin{eqnarray*}
\xi\circ v_1\,&=&\, \xi\circ v_{n+1} + \sum\limits_{i=1}^n p(z_{i},
z_{i+1})\,\xi(h_i)  + \xi\left(\sum\limits_{i=1}^n
d_i\right).\\
&\geq&\, \xi\circ v_{n+1} + \alpha\left(\sum\limits_{i=1}^n p(z_{i},
z_{i+1})\right).
\end{eqnarray*}
From this and using  (B2), we have
\begin{eqnarray*}
\sum\limits_{i=1}^n p(z_{i},
z_{i+1})\,&\leq&\,\frac{1}{\alpha}(\xi\circ v_1-\xi\circ
v_{n+1})\\
&\leq&\,\frac{1}{\alpha}(\xi\circ v_1-\inf \xi\circ
f(S(x_0)))\\
&<&\,+\infty.
\end{eqnarray*}
Thus, $\sum\limits_{i=1}^{\infty}p(z_{i}, z_{i+1})\,<\, +\infty$. By
(B4),  $H$ is $\sigma$-convex, we conclude that
$$\frac{\sum\limits_{i=1}^{\infty} p(z_{i}, z_{i+1})
h_i}{\sum\limits_{j=1}^{\infty} p(z_{j}, z_{j+1})} $$  is
convergent to some  point    $\bar{h}\in H.$
Put
$$h_n^{\prime}:=\frac{\sum\limits_{i=1}^n p(z_{i}, z_{i+1})
h_i}{\sum\limits_{j=1}^n p(z_{j}, z_{j+1})}.$$ Then every
$h_n^{\prime}\in H$ and $h_n^{\prime}\rightarrow\bar{h}$. Remark
that
$$\sum\limits_{i=1}^n p(z_i, z_{i+1}) \geq p(z_1, u) - p(z_{n+1}, u)
> p(z_1, u) -\frac{1}{n+1}.$$
Combining this with (4.4), we have
\begin{eqnarray*}
v_1\,&\in&\, v_{n+1} + \sum\limits_{i=1}^n p(z_i, z_{i+1}) h_i +D\\
&=&\, v_{n+1} +\left(\sum\limits_{i=1}^n p(z_i, z_{i+1})\right)
h_n^{\prime} +D\\
&\subset&\, v_{n+1} + \left(p(z_1, u)-\frac{1}{n+1}\right)
h^{\prime}_n +D\\
&\subset&\, f(z_{n+1}) + \left(p(z_1, u)-\frac{1}{n+1}\right)
h^{\prime}_n +D.
\end{eqnarray*}
By (B3), $f$ is $D$-s.l.m., so $f(z_{n+1})\subset f(u) +D$. Thus,
$$v_1 \in f(u) +\left(p(z_1,u) -\frac{1}{n+1}\right) h^{\prime}_n
+D.\eqno{(4.5)}$$ Since $(p(z_1,u)-1/(n+1)) h^{\prime}_n \rightarrow
p(z_1,u) \bar{h}$ and $f(u) +D$ is closed, by (4.5), we have
$$v_1\in f(u) +p(z_1,u) \bar{h} +D \subset f(u) +p(z_1,u) H +D.$$
Thus, we have shown that $$f(z_1)\subset f(u) + p(z_1,u)H +D,\ \
{\rm i.e.,}\ \ u\in S(z_1) =S(x_{n_0}).$$ This means that assumption
(F) is satisfied.

If the condition that $f$ has $D$-closed values is replaced by one
that  $f$ has $D$-locally closed values, the result remains true.
Since a $\sigma$-convex set $H$ is bounded, the sequence
$(h^{\prime}_n)$ in $H$ is indeed locally convergent to $\bar{h}$.
Thus, the sequence $\left(\left(p(z_1,
u)-\frac{1}{n+1}\right)h^{\prime}_n\right)_{n\in N}$ is locally
convergent to $p(z_1,u)\bar{h}$.  From (4.5), we can deduce that
$v_1\in f(u) + p(z_1, u)\bar{h} +D$ as well since
$f(u)+D$ is locally closed. \hfill\framebox[2mm]{}\\

In Theorem 4.4, if we strengthen the condition that $f$ has
$D$-closed values, then the condition that $H$ is $\sigma$-convex
may be replaced by a weaker one: $H$ is bounded. Let us introduce a
new property on set-valued maps. A set-valued map $f:\,X\rightarrow
2^Y\backslash \{\emptyset\}$ is said to have $(D, H)$-closed values
(or $(D,H)$-locally closed values)  if for any $x\in X$ and any
$\lambda \geq 0$, $f(x) +D+\lambda H$ is closed (or locally
closed). Now, we can give a variant of Theorem 4.4,
which is indeed a generalization of [A pre-order, Theorem 4.2$^{\prime}$].\\

{\bf Theorem 4.5.} \ {\sl Let $(X, {\cal U})$ be a uniform space,
$Y$ be a locally convex space quasi ordered by a convex cone $D$,
$H\subset D\backslash\{0\}$ be a convex set, $f:\, X\rightarrow
2^Y\backslash\{\emptyset\}$ be a set-valued map and $p:\, X\times
X\rightarrow [0, +\infty)$ be a p-distance. Let $x_0\in X$ such that
$S(x_0):=\{x\in X:\, f(x_0)\subset f(x) + p(x_0, x) H +D\} \not=
\emptyset$ and $(X, {\cal U})$ be $S(x_0)$-dynamically complete.
Suppose that the following assumptions are satisfied:

{\rm (B1)} \ $H^{+s}\cap D^+ \not= \emptyset$.

{\rm (B2)} \  $\exists\, \xi\in H^{+s}\cap D^+$ such that $\xi$ is
lower bounded on $f(S(x_0))$.

{\rm (B3)} \  $f$ is $D$-s.l.m.

{\rm (B4$^{\prime}$)} \   $H$ is bounded and $f$ has $(D, H)$-closed
values (or $(D, H)$-locally closed values).

Then, there exists $\hat{x}\in X$ such that

{\rm (a)}  $f(x_0)\subset f(\hat{x}) +p(x_0, \hat{x}) H +D$;

{\rm (b)} $\forall x\in X\backslash\{\hat{x}\},
f(\hat{x})\not\subset f(x) + p(\hat{x}, x) H +D.$}\\

{\bf Proof.} \  As we have seen in the proof of Theorem 4.4, we only
need to prove that assumption (F) in Theorem 3.1$^{\prime}$ is
satisfied.

Let a sequence $(x_n) \subset S(x_0)$ such that $x_n\in S(x_{n-1})$
and
$$\inf\,\xi\circ f(x_n)  < \inf\,\xi\circ f(S(x_{n-1})) +\epsilon_n\
\ \forall n,$$ where every $\epsilon_n >0$ and
$\epsilon_n\rightarrow 0$. As shown in the proof of Theorem 4.4,
there exists $u\in X$ such that $p(x_n, u) \rightarrow 0\
(n\rightarrow\infty)$. Take any fixed $n_0\in N$ and put
$z_1:=x_{n_0}$. As $p(x_k, u) \rightarrow 0\ (k\rightarrow\infty)$,
we may choose a sequence $(z_n)$ from $(x_k)$ such that $p(z_{n+1},
u) < 1/(n+1)$ and $z_{n+1}\in S(z_n),\ \forall n$. Take any $v_1\in
f(z_1)$. As done in the proof of Theorem 4.4, we may choose
$v_{n+1}\in f(z_{n+1}),\ h_n\in H$ and $d_n\in D,\ \forall n$, such
that
$$v_1\,=\, v_{n+1} + \sum\limits_{i=1}^n p(z_{i}, z_{i+1}) h_i
+\sum\limits_{i=1}^n d_i.$$ Thus,
\begin{eqnarray*}
v_1\,&\in&\, f(z_{n+1}) +\left(\sum\limits_{i=1}^n p(z_i,
z_{i+1})\right) H +D\\
&\subset&\, f(z_{n+1}) +\left(p(z_1, u)-p(z_{n+1}, u)\right) H+D\\
&\subset&\, f(z_{n+1}) +\left(p(z_1, u)-\frac{1}{n+1}\right) H +D\\
&\subset&\, f(u) + p(z_1,u) H +D -\frac{1}{n+1} H,\ \ \forall n.
\end{eqnarray*}
Hence, for each $n$, there exists $h_n^{\prime\prime}\in H$ such
that
$$v_1+\frac{1}{n+1} h_n^{\prime\prime}\in f(u) + p(z_1, u) H +D,\ \ \forall
n.$$ By {\rm (B4$^{\prime}$)}, $H$ is bounded, so
$(v_1+\frac{1}{n+1}h_n^{\prime\prime})_{n\in N}$ is (locally)
convergent to $v_1$. Since $f(u) + p(z_1, u) H +D$ is closed
(locally closed), we have
$$v_1\in f(u) +p(z_1,u) H +D.$$
Thus, we have shown that $$f(z_1)\subset f(u) + p(z_1,u)H +D, \  \
{\rm i.e.,}\ \ u\in S(z_1) = S(x_{n_0}).$$\hfill\framebox[2mm]{}\\

{\bf Theorem 4.6.} \ {\sl Let $(X, {\cal U})$ be a uniform space,
$Y$ be a locally convex space whose topology is generated by a
saturated family $\{p_{\alpha}\}_{\alpha\in\Lambda}$ of semi-norms,
$D\subset Y$ be a convex cone, $H\subset D\backslash\{0\}$ be a
convex set, $f:\, X\rightarrow 2^Y\backslash\{\emptyset\}$ be a
set-valued map and $p:\, X\times X\rightarrow [0,+\infty)$ be a
p-distance. Let $x_0\in X$ such that $S(x_0)\not=\emptyset$ and $(X,
{\cal U})$ be $S(x_0)$-dynamically complete. Suppose that the
following assumptions are satisfied:

{\rm (B1$^{\prime}$)} \  $H$ is closed and for each $\alpha\in
\Lambda$, there exists $\xi_{\alpha}\in D^+\backslash\{0\}$ and
$\lambda_{\alpha}>0$ such that $\lambda_{\alpha}p_{\alpha}(h)\leq
\xi_{\alpha}(h),\ \forall h\in H$.

{\rm (B2$^{\prime}$)} \  $f(S(x_0))$ is quasi-bounded from below,
i.e., there exists a bounded set $M$ such that $f(S(x_0))\subset
M+D$.

{\rm (B3)} \  $f$ is $D$-s.l.m.

{\rm (B4$^{\prime\prime}$)} \   $Y$ is $l^{\infty}$-complete {\rm
(see [15, 16])} and $f$ has  $D$-closed values.

Then, there exists $\hat{x}\in X$ such that

{\rm (a)}  $f(x_0)\subset f(\hat{x}) +p(x_0, \hat{x}) H +D$;

{\rm (b)} $\forall x\in X\backslash\{\hat{x}\},
f(\hat{x})\not\subset f(x) + p(\hat{x}, x) H +D.$}\\

{\bf Proof.} \  Since $0\not\in H$ and $H$  is closed, there exists
$\alpha_0\in\Lambda$ and $\eta >0$ such that
$p_{\alpha_0}(h)\geq\eta,\ \forall h\in H$. By {\rm
(B1$^{\prime}$)},
$$\lambda_{\alpha_0}\eta \leq \lambda_{\alpha_0}p_{\alpha_0}(h)\leq
\xi_{\alpha_0}(h),\ \ \forall h\in H.$$ Thus, $\xi_{\alpha_0}\in
D^+\cap H^{+s}$. By {\rm (B2$^{\prime}$)}, $\xi_{\alpha_0}$ is lower
bounded on $f(S(x_0))$. As shown in the proof of Theorem 4.1,
assumptions (D) and (E) in Theorem 3.1$^{\prime}$ are satisfied when
putting $\xi=\xi_{\alpha_0}$. Hence, we only need to prove that
assumption (F) is satisfied for $\xi=\xi_{\alpha_0}$.

Let a sequence $(x_n) \subset S(x_0)$ such that $x_n\in S(x_{n-1})$
and
$$\inf\,\xi_{\alpha_0}\circ f(x_n)  < \inf\,\xi_{\alpha_0}\circ f(S(x_{n-1})) +\epsilon_n\
\ \forall n,$$ where every $\epsilon_n >0$ and
$\epsilon_n\rightarrow 0$. For each $n$, take $y_n\in f(x_n)$ such
that
$$\xi_{\alpha_0}\circ y_n < \inf\,\xi_{\alpha_0}\circ f(S(x_{n-1})) + \epsilon_n.
\eqno{(4.6)}$$ When $m>n$, we have $x_m\in S(x_n)$ and hence
$$y_n \in f(x_n)\subset f(x_m) +p(x_n, x_m) H +D.$$
Thus, there exists $y_{m,n}\in f(x_m),\ h_{m,n}\in H$ and
$d_{m,n}\in D$ such that
$$y_n = y_{m,n} +p(x_n, x_m) h_{m,n} +d_{m,n}.$$
Acting upon the two sides of the above equality by $\xi_{\alpha_0}$,
we have
\begin{eqnarray*}
\xi_{\alpha_0}\circ y_n \,&=&\,\xi_{\alpha_0}\circ y_{m,n} +
p(x_n,x_m)\,
\xi_{\alpha_0}\circ h_{m,n} + \xi_{\alpha_0}\circ  d_{m,n}\\
&\geq&\,\xi_{\alpha_0}\circ y_{m,n} + p(x_n, x_m)\cdot
\lambda_{\alpha_0}\cdot\eta.
\end{eqnarray*}
From this and using (4.6), we have
\begin{eqnarray*}
p(x_n, x_m)\,&\leq&\,\frac{1}{\lambda_{\alpha_0}\cdot
\eta}(\xi_{\alpha_0}\circ y_n -\xi_{\alpha_0}\circ y_{m,n})\\
&\leq&\, \frac{1}{\lambda_{\alpha_0}\cdot \eta}(\xi_{\alpha_0}\circ
y_n -\inf\,\xi_{\alpha_0}\circ
f(S(x_{n-1})))\\
&<&\,\frac{\epsilon_n}{\lambda_{\alpha_0}\cdot \eta}.
\end{eqnarray*}
Thus, $p(x_n, x_m) \rightarrow 0 \ (m>n\rightarrow\infty)$ and
$(x_n)$ is a Cauchy sequence. As shown in the proof of Theorem 4.4,
there exists $u\in X$ such that $p(x_n, u) \rightarrow 0\
(n\rightarrow\infty)$. Take any fixed $n_0\in N$ and put
$z_1:=x_{n_0}$. As $p(x_k, u) \rightarrow 0\ \
(k\rightarrow\infty)$, we may choose a sequence $(z_n)$ from $(x_k)$
such that $p(z_{n+1}, u) < 1/(n+1)$ and $z_{n+1}\in S(z_n),\ \forall
n.$ Take any $v_1\in f(z_1)$. As done in the proof of Theorem 4.4,
we may choose $v_{n+1}\in f(z_{n+1}),\ h_n\in H$ and $d_n\in D,\
\forall n,$ such that
$$v_1 = v_{n+1} + \sum\limits_{i=1}^n p(z_i, z_{i+1}) h_i +
\sum\limits_{i=1}^n d_i.\eqno{(4.7)}$$ For any $\alpha\in \Lambda$,
\begin{eqnarray*}
\xi_{\alpha}\circ v_1\,&\geq&\,\xi_{\alpha}\circ
v_{n+1}+\sum\limits_{i=1}^n p(z_i, z_{i+1})\, \xi_{\alpha}\circ h_i\\
&\geq&\,\xi_{\alpha}\circ v_{n+1}
+\lambda_{\alpha}\,\sum\limits_{i=1}^n p(z_i, z_{i+1})
p_{\alpha}(h_i).
\end{eqnarray*}
From this,
\begin{eqnarray*}
\sum\limits_{i=1}^n p(z_i, z_{i+1})
p_{\alpha}(h_i)\,&\leq&\,\frac{1}{\lambda_{\alpha}}(\xi_{\alpha}\circ
v_1-\xi_{\alpha}\circ v_{n+1})\\
&\leq&\,\frac{1}{\lambda_{\alpha}}(\xi_{\alpha}\circ v_1
-\inf\,\xi_{\alpha}\circ f(S(x_0)))\\
&<&\,+\infty, \hspace{8.5cm}  (4.8)
\end{eqnarray*}
where we use the assumption that $f(S(x_0))$ is quasi-bounded from
below. Thus, $$\sum\limits_{i=1}^{\infty} p(z_i, z_{i+1})\,
p_{\alpha}(h_i) \,<\,+\infty,\ \ \forall \alpha\in \Lambda.$$ By
{\rm (B4$^{\prime\prime}$)}, $Y$ is $l^{\infty}$-complete, so
$\sum\limits_{i=1}^{\infty} p(z_i, z_{i+1}) h_i$ is convergent. On
the other hand, putting $\alpha =\alpha_0$ in (4.8),
$$\eta\cdot \sum\limits_{i=1}^n p(z_i, z_{i+1})\leq \sum\limits_{i=1}^n p(z_i, z_{i+1}) p_{\alpha_0}(h_i)\leq
\frac{1}{\lambda_{\alpha_0}}(\xi_{\alpha_0}\circ v_1
-\inf\,\xi_{\alpha_0}\circ f(S(x_0))) <+\infty.$$ Thus,
$\sum\limits_{i=1}^{\infty} p(z_i, z_{i+1})\,<\,+\infty$. Put
$$\bar{h}:=\, \frac{\sum\limits_{i=1}^{\infty} p(z_{i}, z_{i+1})
h_i}{\sum\limits_{j=1}^{\infty} p(z_{j}, z_{j+1})} \ \ \ {\rm and}\
\ \  h_n^{\prime}:=\frac{\sum\limits_{i=1}^n p(z_{i}, z_{i+1})
h_i}{\sum\limits_{j=1}^n p(z_{j}, z_{j+1})}.$$ Then, every
$h^{\prime}_n\in H$ and $h^{\prime}_n \rightarrow \bar{h}$. Since
$H$ is closed, $\bar{h}\in H$. By (4.7) and (B3), we have
\begin{eqnarray*}
v_1\,&=&\, v_{n+1} + \sum\limits_{i=1}^n p(z_{i}, z_{i+1}) h_i
+\sum\limits_{i=1}^n d_i\\
&=&\,v_{n+1} + \left(\sum\limits_{i=1}^n p(z_{i}, z_{i+1})\right)
h^{\prime}_n
+\sum\limits_{i=1}^n d_i\\
&\in&\, f(z_{n+1}) +\left(\sum\limits_{i=1}^n p(z_i,
z_{i+1})\right) h^{\prime}_n +D\\
&\subset&\, f(u) +\left(\sum\limits_{i=1}^n p(z_i,
z_{i+1})\right) h^{\prime}_n +D\\
&\subset&\, f(u) + (p(z_1, u)-p(z_{n+1}, u))h^{\prime}_n +D\\
&\subset&\, f(u) + \left(p(z_1,u)-\frac{1}{n+1}\right)
h^{\prime}_n+D.
\end{eqnarray*}
Since $(p(z_1, u) -\frac{1}{n+1} )h^{\prime}_n \rightarrow p(z_1, u)
\bar{h}$ and $f(u) +D$ is closed, we have
$$v_1\,\in\, f(u) + p(z_1, u)\bar{h} +D \,\subset\, f(u)+p(z_1, u) H +D.$$
Thus, we have shown that $$f(z_1)\subset f(u) + p(z_1,u)H +D, \  \
{\rm i.e.,}\ \ u\in S(z_1) = S(x_{n_0}).$$\hfill\framebox[2mm]{}\\

Obviously, Theorem 4.6 generalizes [20, Theorem 4.3]. Hence, it is
also a generalization of [13, Theorem 3.5(i)].

We notice that Theorems 4.1 -- 4.3 are set-valued EVPs where
perturbations contain r-distances and Theorems 4.4 -- 4.6 are
set-valued EVPs where perturbations contain p-distances. Finally, we
consider set-valued EVPs where perturbations contain q-distances.
Let us recall the definition  of sequential completeness with
respect to
a q-distance.\\

{\bf Definition 4.4.} \ (see [21]) \  Let $(X, {\cal U})$ be a
uniform space and $p$ be a q-distance on $X$. $(X, {\cal U})$ is
said to be sequentially complete with respect to $p$, if for any
sequence $(x_n)$ in $X$ with $p(x_n, x_m) \rightarrow 0\
(m>n\rightarrow\infty)$, there exists $u\in X$ such that $p(x_n,
u)\rightarrow 0\ (n\rightarrow\infty)$.\\

Moreover, we introduce $S(x)$-dynamical completeness with respect to
a q-distance.\\

{\bf Definition 4.5.} \  As in Definition 4.3, let $(X, {\cal U})$
be a uniform space and $p$ be a q-distance on $X$. Let $S(\cdot):\,
X\rightarrow  2^X\backslash\{\emptyset\}$ be a set-valued map and
$x\in X$. $(X, {\cal U})$ is said to be $S(x)$-dynamically complete
with respect to $p$ if for any sequence $(x_n)\subset S(x)$ such
that $S(x_{n+1})\subset S(x_n)\subset S(x)$ for all $n$ and such
that $p(x_n, x_m) \rightarrow 0\ (m>n\rightarrow\infty)$, there
exists $u\in X$ such that $p(x_n, u) \rightarrow 0 \
(n\rightarrow\infty)$.\\

By checking the proofs of Theorems 4.4 -- 4.6 and using Definition
4.5,  we can easily obtain the following Theorem 4.7, which indeed
includes three theorems corresponding to
Theorems 4.4 -- 4.6.\\

{\bf Theorem 4.7.} \ {\sl In Theorems 4.4, 4.5 and 4.6, if the
assumption that $p$ is a p-distance and $(X, {\cal U})$ is
$S(x_0)$-dynamically complete is replaced by one that $p$ is a
q-distance and $(X, {\cal U})$ is $S(x_0)$-dynamically complete with
respect to $p$, then all the results remain true.}\\

\noindent{\bf References} \vskip 10pt
\begin{description}
\small

\item{[1]} T. Q. Bao, B. S. Mordukhovich, Variational principles for
set-valued mappings with applications to multiobjective
optimization, Control Cybern., 36 (2007), 531-562.

\item{[2]} T. Q. Bao, B. S. Mordukhovich, Relative Pareto minimizers
for multiobjective problems: existence and optimality conditions,
Math. Program, Ser.A, 122 (2010), 301-347.

\item{[3]}  E. M. Bednarczk, D. Zagrodny, Vector variational
principle, Arch. Math. (Basel), 93 (2009), 577-586.

\item{[4]}  F. Flores-Baz\'{a}n, C. Guti\'{e}rrez,  V. Novo,  A
Br\'{e}zis-Browder principle on partially ordered spaces and related
ordering theorems, J. Math. Anal. Appl., 375 (2011), 245-260.

\item{[5]} A. G\"{o}pfert,  C. Tammer and C. Z$\breve{a}$linescu,  On the
vectorial Ekeland's variational principle and minimal point theorems
in product spaces, Nonlinear Anal. 39 (2000),  909-922.

\item{[6]} C. Guti\'{e}rrez, B. Jim\'{e}nez, V. Novo, A set-valued
Ekeland's variational principle in vector  optimization, SIAM J.
Control. Optim., 47 (2008), 883-903.

\item{[7]} T. X. D. Ha, Some variants of the Ekeland  variational
principle for a set-valued map, J. Optim. Theory Appl., 124 (2005),
187-206.

\item{[8]}  A. H. Hamel, Equivalents to Ekeland's variational
principle in uniform spaces, Nonlinear Anal. 62 (2005), 913-924.

\item{[9]} O. Kada, T. Suzuki, W. Takahashi, Nonconvex minimization
theorems and fixed point theorems in complete metric spaces, Math.
Japon., 44 (1996), 381-391.

\item{[10]}   P. Q. Khanh, D. N. Quy, On generalized Ekeland's variational
principle and equivalent formulations for set-valued mappings, J.
Glob. Optim., 49 (2011), 381-396.

\item{[11]}  P. Q. Khanh, D. N. Quy, Versions of Ekeland's
variational principle involving set perturbations, J. Glob. Optim.,
57 (2013), 951-968.

\item{[12]}  L. J. Lin, W. S. Du, Ekeland's variational principle,
minimax theorems and existence of nonconvex equilibria in complete
metric spaces, J. Math. Anal. Appl., 323 (2006), 360-370.

\item{[13]}  C. G. Liu, K. F. Ng, Ekeland's variational principle for
set-valued functions, SIAM J. Optim., 21 (2011), 41-56.

\item{[14]}  A. B. N\'{e}meth, A nonconvex vector minimization problem,
Nonlinear Anal. 10 (1986),  669-678.

\item{[15]}  P. P\'{e}rez Carreras, J. Bonet, Barrelled Locally Convex Spaces,
North-Holland, Amsterdam, 1987.

\item{[16]} J. H. Qiu, Local completeness and dual local quasi-completeness,
Proc. Amer. Math. Soc. 129 (2001) 1419-1425.

\item{[17]} J. H. Qiu, On Ha's version of set-valued Ekeland's
variational principle, Acta Math. Sinica, English Series, 28 (2012),
717-726.

\item{[18]} J. H. Qiu, Set-valued quasi-metrics and a general
Ekeland's variational principle in vector optimization, SIAM J.
Control  Optim., 51 (2013), 1350-1371.

\item{[19]} J. H. Qiu, The domination property for efficiency and
Bishop-Phelps theorem in locally convex spaces, J. Math. Anal.
Appl., 402 (2013), 133-146.

\item{[20]} J. H. Qiu, A pre-order principle and set-valued Ekeland
variational principle, arXiv: 1311.4951[math.FA].

\item{[21]} J. H. Qiu, F. He, p-distances, q-distances and a
generalized Ekeland's variational principle in uniform spaces,
 Acta Math. Sinica, English Ser., 28  (2012), 235-254.

\item{[22]} C. Tammer, C. Z$\breve{a}$linescu, Vector variational
principle for set-valued functions, Optimization, 60 (2011),
839-857.

\item{[23]} C. Z$\breve{a}$linescu, Convex Analysis  in   General
Vector Spaces, World Sci., Singapore, 2002.

\end{description}
\end{document}